\documentclass{ACmod}

\DeclareFontFamily{U}{rsf}{}
\DeclareFontShape{U}{rsf}{m}{n}{
  <5> <6> rsfs5 <7> <8> <9> rsfs7 <10-> rsfs10}{}
\DeclareMathAlphabet{\mathscr}{U}{rsf}{m}{n}

\makeatletter
\def\operator@font{\sf}
\makeatother

\newcommand{\cO}{{\mathscr O}}

\newcommand{\HKR}{{\mathsf{HKR}}}

\DeclareMathOperator{\Sym}{Sym}

\DeclareMathOperator{\HH}{HH}
\DeclareMathOperator{\HT}{HT}

\DeclareMathOperator{\Td}{td}

\DeclareMathOperator{\id}{id}

\DeclareMathOperator{\Ext}{Ext}

\DeclareMathOperator{\Hom}{Hom}

\renewcommand{\phi}{\varphi}

\usepackage{tikz-cd}
\usepackage{amsmath}

\usepackage{xypic}
\usepackage{tikz}
\usepackage[tocflat]{tocstyle}
\usetikzlibrary{positioning}
\usepackage{amscd}
\usetikzlibrary{arrows,shapes,positioning}
\usetikzlibrary{decorations.markings}
\tikzstyle arrowstyle=[scale=1]
\tikzstyle directed=[postaction={decorate,decoration={markings,
    mark=at position .65 with {\arrow[arrowstyle]{stealth}}}}]
\tikzstyle reverse directed=[postaction={decorate,decoration={markings,
    mark=at position .65 with {\arrowreversed[arrowstyle]{stealth};}}}]
\usepackage[tocflat]{tocstyle}
\tikzset{cross/.style={cross out, draw=black, minimum size=2*(#1-\pgflinewidth), inner sep=0pt, outer sep=0pt},
cross/.default={1pt}}

\begin{document}

\title{When are two HKR isomorphisms equal?}

\author[Shengyuan Huang]{%
Shengyuan Huang}

\address{
\textsc{Shengyuan Huang: School of Mathematics,
University of Birmingham
\newline
Edgbaston, Birmingham, B15 2TT, UK}
\newline
Email address: \texttt{s.huang.5@bham.ac.uk}}

\begin{abstract}
{\sc Abstract:}
 Let $X\hookrightarrow S$ be a closed embedding of smooth schemes which splits to first order. An HKR isomorphism is an isomorphism between the shifted normal bundle $\mathbb{N}_{X/S}[-1]$ and the derived self-intersection $X\times^R_SX$. Given two different first order splittings of a closed embedding, one can obtain two HKR isomorphisms using a construction of Arinkin and C\u{a}ld\u{a}raru. A priori, it is not known if the two isomorphisms are equal or not. We define the generalized Atiyah class of a vector bundle on $X$ associated to a closed embedding and two first order splittings. We use the generalized Atiyah class to give sufficient and necessary conditions for when the two HKR isomorphisms are equal over $X$ and over $X\times X$ respectively. When $i$ is the diagonal embedding, there are two natural projections from $X\times X$ to $X$. We show that the HKR isomorphisms defined by the two projections are equal over $X$, but not equal over $X\times X$ in general. 
 

\end{abstract}

\maketitle
\setcounter{tocdepth}{1}
\tableofcontents
\vfill

\section{Introduction}
\paragraph Let $X$ be a smooth algebraic variety over a field of characteristic zero. There is an HKR isomorphism~\cite{S,Y} in the derived category of $X$
$$\Delta^*\Delta_*\cO_X\cong \Sym_{\cO_X}(\Omega_X[1]),$$
where $\Delta: X\hookrightarrow X\times X$ is the diagonal embedding.

One can interpret the isomorphism above in terms of an isomorphism of derived schemes. The structure complex of the derived self-intersection $X\times^R_{X\times X}X$ is $\Delta^*\Delta_*\cO_X$. The structure complex of the shifted tangent bundle $\mathbb{T}_{X}[-1]$ is $\Sym_{\cO_X}(\Omega_X[1])$. The HKR isomorphism above can be viewed as an isomorphism between the shifted tangent bundle and the derived self-intersection
$$\mathbb{T}_{X}[-1]\cong X\times^R_{X\times X}X.$$

\paragraph One can replace the diagonal embedding by an arbitrary closed embedding $i:X\hookrightarrow S$ of smooth schemes and consider the derived self-intersection $X\times^R_SX$. The embedding $i$ factors as
$$
\xymatrix{
X\ar@{^{(}->}[r]^{\mu} & X^{(1)}_S\ar@{^{(}->}[r]^{\nu} & S.
}
$$
where $X^{(1)}_S$ is the first order neighborhood of $X$ in $S$.  We say $i$ splits to first order if and only if the map $\mu$ is split, i.e., there exists a map of schemes $\varphi: X^{(1)}_S\rightarrow X$ such that $\varphi\circ\mu=id$. 

\paragraph\label{par 1.1} There is a bijection between first order splittings of $i$ and splittings of the short exact sequence below \cite[20.5.12 (iv)]{G}
$$
\xymatrix{
 0\ar[r] & T_X \ar[r] & T_S|_X \ar[r] & N_{X/S}\ar[r]\ar@{-->}@/_1.0pc/[l] & 0.
 }
$$
The bijection above is canonical, so we use the same notation $\varphi$ for the splitting of the bundles and the splitting map $X^{(1)}_S\rightarrow X$ of schemes.

We assume that $i$ splits to first order throughout this paper. Let $\mathbb{N}_{X/S}[-1]$ be the total space of the shifted normal bundle $N_{X/S}[-1]$. Namely, the space $\mathbb{N}_{X/S}[-1]$ is a derived scheme whose structure complex is $\Sym_{\cO_X}(N^{\vee}_{X/S}[1])$. From a fixed first order splitting, Arinkin and C\u{a}ld\u{a}raru~\cite{AC} constructed an isomorphism
$$\HKR_{\varphi}:\mathbb{N}_{X/S}[-1]\cong X\times^R_SX$$
between the shifted normal bundle $\mathbb{N}_{X/S}[-1]$ and the derived self-intersection $X\times^R_SX$. When $i$ is the diagonal embedding, the normal bundle $N_{X/S}$ is isomorphic to the tangent bundle $T_X$.

\paragraph{\bf A positive answer to Arinkin-C\u{a}ld\u{a}raru's question.} When $i$ is the diagonal embedding $\Delta:X\hookrightarrow X\times X$, there are two first order splitting $\pi_1$ and $\pi_2$ obtained from the two projections $p_1$ and $p_2:X\times X\rightarrow X$. In this case, Yekutieli~\cite{Y} defined a complete bar resolution to compute the HKR isomorphism for a general scheme $X$. It has been shown that the HKR isomorphism obtained by Yekutieli is equal to the HKR isomorphism defined by the second projection~\cite{C}.

Arinkin and C\u{a}ld\u{a}raru~\cite[Paragraphs 1.21-1.24]{AC} asked if the two HKR isomorphisms $\HKR_{\pi_1}$ and $\HKR_{\pi_2}$ are equal over $X$. Grivaux provided~\cite{Gri'} a positive answer to this question. Corollary~\ref{cor over X diagonal} provides a different proof of Grivaux's result.

\paragraph{\bf Grivaux's question.} In the case of general embedding $X\hookrightarrow S$, two different first order splittings $\varphi_1$ and $\varphi_2$ define two HKR isomorphisms $\HKR_{\varphi_1}$ and $\HKR_{\varphi_2}:$ $i^*i_*\cO_X \cong \Sym(N^\vee_{X/S}[1])$. The composite map $\HKR_{\varphi_1}\circ\HKR^{-1}_{\varphi_2}$ defines an automorphism of $\Sym(N^\vee_{X/S}[1])$. Grivaux asked if we can compute this automorphism explicitly and answered this question when the codimension of $X$ in $Y$ is two~\cite[Theorem 4.17]{Gri'}. Theorem A gives a complete answer to this question.

\paragraph Before we state the main theorems, we need to explain a technical detail. One can consider the HKR isomorphism over $X$ or over $X\times X$. An HKR isomorphism $\mathbb{N}_{X/S}[-1]\cong X\times^R_SX$ over $X$ is equivalent to an algebra isomorphism $i^*i_*(\cO_X)\cong\Sym (N^{\vee}_{X/S}[1])$ of the structure complexes of the two derived schemes. An isomorphism $\mathbb{N}_{X/S}[-1]\cong X\times^R_SX$ over $X\times X$ is more complicated. In this case, we can view the structure complex of $\mathbb{N}_{X/S}[-1]$ as an $\cO_{X\times X}$-module, i.e., an $\cO_{X}$-bimodule. This bimodule is a kernel that represents the dg functor $i^*i_*(-): \mathsf{D}(X) \rightarrow \mathsf{D}(X)$ from the dg enhancement $\mathsf{D}(X)$ of the derived category of $X$ to itself. It turns out that an isomorphism over $X\times X$ is equivalent to an isomorphism $i^*i_*(-)\cong\Sym (N^{\vee}_{X/S}[1])\otimes(-)$ of dg functors~\cite{ACH}. An isomorphism over $X\times X$ implies that the isomorphism is also over $X$. More details about the two natural base schemes have been explained in~\cite{ACH}. From a fixed first order splitting, the HKR isomorphism constructed in~\cite{AC} is over $X\times X$ as explained in Section 2.

\paragraph
One can obtain two HKR isomorphisms from two different first order splittings. We construct a cohomology class below associated to the two splittings and we give sufficient and necessary conditions for the two HKR isomorphisms to be equal over $X$ and over $X\times X$ respectively.

\begin{Definition} \label{seq}
Let $\varphi_1$ and $\varphi_2$ be two first order splittings of the map $i:X\hookrightarrow S$. We construct a cohomology class $\alpha_{\varphi_1,\varphi_2}(E)$ associated to $\varphi_1$, $\varphi_2$, and a vector bundle $E$ on $X$ in the following way. 

The difference $\varphi_1-\varphi_2$ is a map $\Omega_X\rightarrow N^\vee_{X/S}$ due to the bijiection in Paragraph~\ref{par 1.1}. The class $\alpha_{\varphi_1,\varphi_2}(E)$ is defined as the following composite map
$$E\rightarrow E\otimes\Omega_X[1]\rightarrow E\otimes N^\vee_{X/S}[1],$$
where the second map $E\otimes\Omega_X[1]\rightarrow E\otimes N^\vee_{X/S}[1]$ is $\id_E\otimes(\varphi_1-\varphi_2)$, the first map is the Atiyah class of $E$. We call the class above the {\emph generalized Atiyah class} of $E$ associated to $\varphi_1$ and $\varphi_2$. 
\end{Definition}

When $i$ is the diagonal embedding, and $\varphi_1$ and $\varphi_2$ are the first order splittings corresponding to the projections onto the first and second factor, one can check that $N^\vee_{X/S}$ is isomorphic to the cotangent bundle $\Omega_{X}$ and $\varphi_1-\varphi_2$ can be identified as the identity map from $\Omega_X$ to itself. The class is the classical Atiyah class $\mathsf{at}(E)$ of $E$.\\


Before we state the main theorems we need a few classes below. The class $\alpha_{\varphi_1,\varphi_2}(N^\vee_{X/S})$ can be viewed as a map $N^\vee_{X/S}\rightarrow {N^\vee_{X/S}}^{\otimes 2}[1]$. Denote the composite map
$$N^\vee_{X/S}\stackrel{\alpha_{\varphi_1,\varphi_2}}{\longrightarrow} {N^\vee_{X/S}}^{\otimes 2}[1]\stackrel{\cdot\frac{1}{2}}{\rightarrow}{N^\vee_{X/S}}^{\otimes 2}[1]\rightarrow\wedge^2{N^\vee_{X/S}}[1]$$
by $\alpha^{\mathsf{antisym}}_{\varphi_1,\varphi_2}(N^\vee_{X/S})$, where the map in the middle is the multiplication by $\frac{1}{2}$, and the last map is the natural projection.

For a first order splitting $\varphi_1$, we define a class $\Phi^k_{\varphi_1}(E):\mu_* E \rightarrow\mu_* ({N^\vee_{X/S}}^{\otimes k}\otimes E)[k]$ for all $k\geq 1$ in Paragraph~\ref{Atiyah}. Due to the complexity of the definition, we do not define the class here and readers can find the definition of the class later in Paragraph~\ref{Atiyah}. We can pushforward the class by $\varphi_2$, so we get $\varphi_{2*}\Phi^k_{\varphi_1}(E):E \rightarrow ({N^\vee_{X/S}}^{\otimes k}\otimes E)[k]$. When $k$ is equal to one, this class is nothing but $\alpha_{\varphi_1,\varphi_2}(E)$.
\paragraph{\bf Theorem A.}
{\em 
Let $i:X\hookrightarrow S$ be a closed embedding of schemes and assume $i$ splits to first order. Let $\varphi_1$ and $\varphi_2$ be two first order splittings. The two different first order splittings define two HKR isomorphisms $\HKR_{\varphi_1}$ and $\HKR_{\varphi_2}:$ $i^*i_*\cO_X \cong \Sym(N^\vee_{X/S}[1])$. The composite map $\HKR_{\varphi_1}\circ\HKR^{-1}_{\varphi_2}$ defines an automorphism of $\Sym(N^\vee_{X/S}[1])$. Write the automorphism in the form of a matrix below. Then the matrix is unipotent
$$\mbox{\small
\xymatrix{
  & N^\vee_{X/S} & \wedge^2 N^\vee_{X/S}[1] & \wedge^3 N^\vee_{X/S}[2] & \cdots\\
N^\vee_{X/S} &  \id &  \alpha^{\mathsf{antisym}}_{\varphi_1,\varphi_2}(E)                &                   \cdots        &    \cdots   \\
\wedge^2 N^\vee_{X/S}[1] & 0  &   \id                  &     \cdots                &  \cdots     \\
\wedge^3 {N^\vee_{X/S}}[2] & 0 & 0 & \id &\\
\cdots & & & &,}
}
$$
where the maps on the diagonal are the identity maps. The $(p,p+k)$-th entry
$$  \wedge^{p-1}{N^\vee_{X/S}}\rightarrow  \wedge^{p+k-1}{N^\vee_{X/S}}[k]$$
is the composite map
$$ \wedge^{p-1}{N^\vee_{X/S}}\hookrightarrow  {N^\vee_{X/S}}^{\otimes p-1}$$
$$\xymatrix{ \ar[rrr]^{\varphi_{2*}\Phi^k_{\varphi_1}( {N^\vee_{X/S}}^{\otimes p-1})~~~~~~~~} &  &  &  {N^\vee_{X/S}}^{\otimes k+p-1}[k]
}$$
$$\stackrel{\frac{1}{(k+p-1)!}}{\longrightarrow} {N^\vee_{X/S}}^{\otimes k+p-1}[k]\rightarrow\wedge^{p+k-1}{N^\vee_{X/S}}[k],$$
where the third arrow in the composite map above is the multiplication by $\frac{1}{(k+p-1)!}$.
}

\paragraph{\bf Theorem B.}
{\em We are in the same setting of Theorem A. 

\begin{enumerate}
\item The two HKR isomorphisms $\HKR_{\varphi_1}$ and $\HKR_{\varphi_2}$ are equal over $X\times X$ if and only if the class $\alpha_{\varphi_1,\varphi_2}(E)$ vanishes for all $E$. 

\item The two isomorphisms are equal over $X$ if and only if the class $\alpha^{\mathsf{antisym}}_{\varphi_1,\varphi_2}(N^\vee_{X/S})$ vanishes.
\end{enumerate}
}

When the map $i$ is the diagonal embedding, the two splittings are $\pi_1$ and $\pi_2$, and $N^\vee_{X/S}$ is the cotangent bundle $\Omega_X$, the class $\alpha_{\pi_1,\pi_2}(N^\vee_{X/S})$ is the Atiyah class of the cotangent bundle. Note that the Atiyah class of the cotangent bundle is symmetric and the class $\alpha_{\pi_1,\pi_2}^{\mathsf{antisym}}(N^\vee_{X/S})$ is the anti-symmetric part of the Atiyah class. Therefore the class $\alpha_{\pi_1,\pi_2}^{\mathsf{antisym}}(N^\vee_{X/S})$ always vanishes in the case of diagonal embedding. We can conclude that the two HKR isomorphisms defined by the two projections $\pi_1$ and $\pi_2$ are equal over $X$ which is Corollary~\ref{cor over X diagonal}. In the case of general embedding $X\hookrightarrow S$, it is not known whether the class $\alpha_{\varphi_1,\varphi_2}(N^\vee_{X/S})$ is symmetric or not.

\paragraph The author would like to point out that Arinkin, C\u{a}ld\u{a}raru, and Hablicsek~\cite{ACH} provided another way to construct an HKR isomorphism from a fixed first order splitting. This construction was also obtained by Grivaux~\cite{Gri} in the context of differential geometry. From a fixed first order splitting, it is not known whether the isomorphisms obtained from the two different constructions in~\cite{AC} and~\cite{ACH} are equal or not. It is highly possible that the two isomorphisms obtained from the two different constructions are equal, but we only consider the first construction in~\cite{AC} throughout this paper.

\paragraph{\bf Plan of the paper.}

In Section 2 we recall the construction of the HKR isomorphisms from a first order splitting.

In Section 3 we provide an alternative definition of the generalized Atiyah class. Then we study the properties of the generalized Atiyah class.

In Section 4, for $k\geq1$, we define two collection of classes $\Phi^k_{\varphi}(E)$ and $\Psi^k_{\varphi}(E)$ from a vector bundle $E$ and a first order splitting $\varphi$. We relate the classes with the construction in Section 2 and then prove the first part of Theorem B.

In Section 5 we prove Theorem A and the second part of Theorem B follows from Theorem A immediately.

\paragraph {\bf Conventions.} All the schemes in this paper are smooth over a field of characteristic zero.

\paragraph {\bf Acknowledgments.} The author is grateful to Kai Xu for helpful discussions and for proving Proposition~\ref{prop universal}. We would like to thank Tyler Kelly and the anonymous referee for their useful suggestions on writing.

The author was partially supported by the UKRI Future Leaders Fellowship through grant number MR/T01783X/1.

\section{The background on HKR isomorphisms}

In this section we recall the constructions of HKR isomorphism $\HKR_{\varphi}$ in~\cite{AC} from a chosen first order splitting $\varphi$. In the construction one build an explicit resolution of $\mu_*\cO_X$ as an $\cO_{X^{(1)}_S}$-algebra using $\varphi$. The explicit resolution is crucial and will be used throughout this paper.

\paragraph\label{resolution} Fix a first order splitting $\varphi$. We recall the construction of HKR isomorphism $\mathbb{N}_{X/S}[-1]\cong X\times^R_SX$ in~\cite{AC} from the given $\varphi$.

To define an isomorphism
$\mathbb{N}_{X/S}[-1]\cong X\times^R_SX$ over $X$, it suffices to define an algebra isomorphism $$i^*i_*\cO_X\cong\Sym(N^\vee_{X/S}[1])$$

on the structure complexes of both derived schemes. The map is defined as the composite map
$$\mbox{\small
\xymatrix{
\mu^*\nu^*\nu_*\mu_*\cO_X\ar[r] & \mu^*\mu_*\cO_X\ar[r]^{\cong} & \mathrm{T}^c(N^\vee_{X/S}[1])\ar[r]^{\exp} & \mathrm{T}(N^\vee_{X/S}[1])\ar[r] &\Sym(N^\vee_{X/S}[1]).
}
}
$$
The leftmost map is given by the counit of the adjunction $\nu^*\dashv\nu_*$.  The map $\exp$ is multiplying by $1/k!$ on the degree $k$ piece, and the last one is the natural projection map.  The $\mathrm{T}^c(N^\vee_{X/S}[1])$ is the free coalgebra on $N^\vee_{X/S}[1]$ with the shuffle product structure, and $\mathrm{T}(N^\vee_{X/S}[1])$ is the tensor algebra on $N^\vee_{X/S}[1]$.  The isomorphism $\mu^*\mu_*\cO_X\cong\mathrm{T}^c(N^\vee_{X/S}[1])$ in the middle is nontrivial and needs more explanation.  With the splitting $\varphi$ one can build an explicit resolution of $\mu_*\cO_X$ as an $\cO_{X^{(1)}_S}$-algebra
$$
 \xymatrix{
 (\mathrm{T}^c(\varphi^*N^\vee_{X/S}[1]),d)\ar[r] & \mu_*\cO_{X},
 }
 $$
where $(\mathrm{T}^c(\varphi^*N^\vee_{X/S}[1]),d)$ is the free coalgebra on $\varphi^*N^\vee_{X/S}[1]$ with the shuffle product structure and a differential $d$.  The construction of the resolution is as follows.

Consider the short exact sequence

$$0\rightarrow\mu_*N^\vee_{X/S}\rightarrow \cO_{X^{(1)}_S}\rightarrow\mu_*\cO_X\rightarrow0.$$

For a vector bundle $E$, tensor the sequence with $\varphi_1^*(E)$. We get
$$0\rightarrow\mu_*(N^\vee_{X/S}\otimes E)\rightarrow \varphi_1^*(E)\rightarrow\mu_* E\rightarrow0 $$
due to the projection formula.

Taking $E$ to be $(N^{\vee}_{X/S})^{\otimes k}$ for all nonnegative integers $k$, we get a family of short exact sequences

$$0\rightarrow\mu_*({N^\vee_{X/S}}^{\otimes k+1})\rightarrow \varphi^*({N^\vee_{X/S}}^{\otimes k})\rightarrow\mu_* ({N^\vee_{X/S}}^{\otimes k})\rightarrow0 $$
for all $k$. Stringing together these exact sequences for successive values of $k$, we get the desired resolution $(\mathrm{T}^c(\varphi^*N^\vee_{X/S}[1]),d)$ of $\mu_*\cO_X$

$$\cdots\rightarrow \varphi^*({N^\vee_{X/S}}^{\otimes k+1})\rightarrow \varphi^*({N^\vee_{X/S}}^{\otimes k})\rightarrow\cdots\rightarrow\cO_{X^{(1)}_S}\rightarrow\mu_*\cO_X\rightarrow 0.$$

The differential vanishes once we pull this resolution back on $X$ via $\mu$, so we get the desired isomorphism $\mu^*\mu_*\cO_X\cong\mathrm{T}^c(N^\vee_{X/S}[1])$.

For any sheaf $E$ on $X$, we tensor the resolution above by $\varphi^*E$.  Using the projection formula and $\varphi\circ\mu=id$, one can show that we get a resolution of $\mu_*E$
$$
(\mathrm{T}^c(\varphi^*N^\vee_{X/S}[1])\otimes\varphi^*E,d)\rightarrow\mu_*E.
$$
Denote the resolution above by $T_{E,\varphi}$. The same argument shows that $i^*i_*(E)\cong E\otimes\Sym(N^\vee_{X/S}[1])$, i.e., that $i^*i_*(-)\cong(-)\otimes\Sym(N^\vee_{X/S}[1])$ as dg functors. All the constructions above are canonical except for the isomorphism $I_{\varphi}(E):\mu^*\mu_*(E)\cong \mathrm{T}(N^\vee_{X/S}[1])\otimes E$ which depends on the splitting.

The author would like to address that $\mathrm{T}^c(N^\vee_{X/S}[1])$ is equal to $\mathrm{T}(N^\vee_{X/S}[1])$ as chain complexes, but we would like to remember the commutative algebra structure on $\mathrm{T}^c(N^\vee_{X/S}[1])$ and the isomorphism $\mu^*\mu_*(\cO_X)\cong\mathrm{T}^c(N^\vee_{X/S}[1])$ is an isomorphism of algebras. We omit the superscript $c$ for the isomorphism $\mu^*\mu_*(E)\cong\mathrm{T}(N^\vee_{X/S}[1])\otimes E$ because there are no algebra structures on both sides for general $E$.

\section{The generalized Atiyah class}

In this section we give an equivalent definition of the generalized Atiyah class and study its properties.

\begin{Proposition}\label{prop universal}
Let $\varphi_1$ and $\varphi_2$ be two first order splittings of the map $i:X\hookrightarrow S$. We construct a cohomology class associated to $\varphi_1$, $\varphi_2$, and a vector bundle $E$ on $X$ in the following way. Consider the short exact sequence

$$0\rightarrow\mu_*N^\vee_{X/S}\rightarrow \cO_{X^{(1)}_S}\rightarrow\mu_*\cO_X\rightarrow0.$$

Tensor the sequence with $\varphi_1^*(E)$. We get
$$0\rightarrow\mu_*(N^\vee_{X/S}\otimes E)\rightarrow \varphi_1^*(E)\rightarrow\mu_* E\rightarrow0 $$
due to the projection formula. Then we pushforward the exact sequence by $\varphi_2$

$$0\rightarrow N^\vee_{X/S}\otimes E\rightarrow \varphi_{2*}\varphi_1^*(E)\rightarrow E\rightarrow0. $$

The sequence above defines an extension class in $\Ext^1(E,E\otimes N^\vee_{X/S})$. This class is equal to the generalized Atiyah class $\alpha_{\varphi_1,\varphi_2}(E)$.

\end{Proposition}

\begin{Proof}
The two first order splittings defines a map $\Sigma=(\varphi_1,\varphi_2):X^{(1)}_S\rightarrow X\times X$. Because of the map $\Sigma$, it induces a map of normal bundles
$$N_{X/S}=N_{X/X^{(1)}_S}\rightarrow N_{X/X\times X}=T_X.$$
One can check that the map is dual to $\varphi_1-\varphi_2$.
Let $p_1$ and $p_2$ be the two projections $X\times X\rightarrow X$. The diagram
$$
\xymatrix{
X^{(1)}_S \ar[rr]^{\Sigma}\ar[dr]^{\varphi_1}\ar[drrr]^{\varphi_2} &  & X\times X\ar[dr]^{p_2}\ar[dl]^{p_1} &\\
  & X & & X}
$$

is commutative due to the definition of $\Sigma$. For any vector bundle $E$, there is a morphism of short exact sequences
$$
\xymatrix{
0 \ar[r] & E\otimes\Omega_X\ar[r]\ar[d]_{\id\otimes(\varphi_1-\varphi_2)} & p_{2*}p^*_1 E\ar[r]\ar[d] & E\ar[r]\ar[d]_{\id} & 0\\
0 \ar[r] & E\otimes N^\vee_{X/S}\ar[r] &  \varphi_{2*}\varphi^*_1 E\ar[r] & E\ar[r] & 0,
}
$$
where the vertical arrow in the middle is defined by adjunction and the equalities $\varphi_1=p_1\circ\Sigma$ and $\varphi_2=p_2\circ\Sigma$. The top line is the exact sequence that defines the Atiyah class of $E$ and the bottom line is the exact sequence in this proposition.
\end{Proof}

Due to the proposition above and \cite{CG}[Paragraph 4.1.3], one can conclude that $\varphi_1^*E$ is isomorphic to $\varphi_2^*E$ if and only if the class $\alpha_{\varphi_1,\varphi_2}(E)$ vanishes. The proposition above implies an interesting result below.
\begin{Corollary}
Let $E$ be a vector bundle on $X$ whose Atiyah class is zero. Then for any closed embedding $i:X\hookrightarrow S$ and any two first order splittings $\varphi_1$ and $\varphi_2$, the two bundles $\varphi_1^*E$ and $\varphi_2^*E$ are isomorphic.
\end{Corollary}

\paragraph \label{num}
The set of all first order splittings
$$
\xymatrix{
 0\ar[r] & T_X \ar[r] & T_S|_X \ar[r] & N_{X/S}\ar[r]\ar@{-->}@/_1.0pc/[l] & 0.
 }
$$
is a $\Hom_{\cO_X}(N_{X/S},T_X)=\Hom_{\cO_X}(\Omega_X, N^\vee_{X/S})$-torsor. We can identify the set of splittings with $\Hom_{\cO_X}(\Omega_X, N^\vee_{X/S})$ by choosing a first order splitting $\varphi$. Then we get a map
$$\Theta_{\varphi}:\Hom_{\cO_X}(\Omega_X, N^\vee_{X/S})\rightarrow \Ext^1(E,E\otimes N^\vee_{X/S}) $$
as follows. For any element $x\in\Hom_{\cO_X}(\Omega_X, N^\vee_{X/S})$, we get another splitting $\varphi+x$. The element $x$ is mapped to the class $\alpha_{\varphi+x,\varphi}(E)\in\Ext^1(E,E\otimes N^\vee_{X/S})$. It is very natural to expect that the map $\Theta_\varphi$ is a linear map between vector spaces and that it does not depend on $\varphi$. We prove this statement in Proposition~\ref{linear}.

\begin{Proposition}\label{linear}
Fix a first order splitting $\varphi$ and it identifies the set of splittings with $\Hom(\Omega_X,N^\vee_{X/S})$. The map $\Theta_{\varphi}$ is a linear map between vector spaces.
\end{Proposition}
\begin{Proof}
Proposition~\ref{prop universal} shows that $\Theta_\varphi(x)=(\id_E\otimes x)\circ\mathsf{at}(E)$. It is clear that this map is linear and it does not depend on the splitting $\varphi$ we choose.
\end{Proof}

Proposition~\ref{linear} has an interesting application in the case of diagonal embedding. Fix a first order splitting $\varphi=\pi_1$, and then the set of first order splittings is identified with $\Hom(\Omega_X,\Omega_X)$. One can show that the difference of the two first order splittings $\pi_1$ and $\pi_2$ is the identity map $\id_{\Omega_X}$. Similarly $\frac{1}{2}\id_{\Omega_X}$ corresponds to a new splitting $\frac{\pi_1+\pi_2}{2}$
$$
\xymatrix{
 0\ar[r] & \Omega_X \ar[r] & \Omega_X\oplus\Omega_X \ar[r]\ar@{-->}@/_1.0pc/[l]_{\frac{\pi_1+\pi_2}{2}} &\Omega_X \ar[r] & 0.
 }
$$

The map $\Theta_{\varphi}$ sends $\id_{\Omega_X}$ to the Atiyah class $\mathsf{at}(E)$.
The linearity of the map implies that $\Theta_{\varphi}(\frac{1}{2}\id_{\Omega_X})=\frac{1}{2}\mathsf{at}(E)$, i.e., the equality
$\alpha_{\pi_1,\frac{\pi_1+\pi_2}{2}}(E)=\frac{1}{2}\mathsf{at}(E)$.

\begin{Proposition}\label{prop derivation}
Let $E$ and $F$ be two vector bundles on $X$. The class $\alpha_{\varphi_1,\varphi_2}$ satisfies the equality
$$\alpha_{\varphi_1,\varphi_2}(E\otimes F)=\id_E\otimes\alpha_{\varphi_1,\varphi_2}(F)+\alpha_{\varphi_1,\varphi_2}(E)\otimes\id_F.$$
\end{Proposition}

\begin{Proof}
The Atiyah class satisfies~\cite{M}
$$\mathsf{at}(E\otimes F)=\id_E\otimes\mathsf{at}(F)+\mathsf{at}(E)\otimes\id_F.$$ Then use Proposition~\ref{prop universal}.
\end{Proof}

\section{Proof of Theorem B (1)}
In this section we define cohomology classes $\Phi^k_{\varphi}(E)$ and $\Psi^K_{\varphi}(E)$ for a first order splitting $\varphi$ and a vector bundle $E$. We explain that the classes are related to the explicit resolution of $\mu_*\cO_{X}$ in Section 2. We study the properties of the classes and then use the properties to prove the first part of Theorem B.

\paragraph The exact sequence $$0\rightarrow\mu_*(N^\vee_{X/S}\otimes E)\rightarrow \varphi^*(E)\rightarrow\mu_* E\rightarrow0 $$ is crucial in this paper.
It defines a map $\Phi^1_{\varphi}(E):\mu_* E\rightarrow \mu_*(N^\vee_{X/S}\otimes E)[1]$. Pushing forward the class onto $S$, we get a map $\Psi^1_{\varphi}(E):i_* E\rightarrow i_*(N^\vee_{X/S}\otimes E)[1]$. We define a collection of maps $\Phi^k_{\varphi}(E)$ and $\Psi^k_{\varphi}(E)$ in Paragraph~\ref{Atiyah}.

Given two splittings $\varphi_1$ and $\varphi_2$, one sees that $\varphi_{2*}\Phi^1_{\varphi_1}(E)$ is equal to $\alpha_{\varphi_1,\varphi_2}(E)$ defined in Definition~\ref{seq} due to Proposition~\ref{prop universal}.

When $i$ is the diagonal embedding $\Delta:X\hookrightarrow X\times X$, $E$ is the structure sheaf $\cO_X$, and $\varphi=\pi_2$ is the first order splitting obtained from the projection $p_2:X\times X \rightarrow X$ onto the second factor. The class $\Psi^1_{\pi_2}(\cO_X)$ is called the universal Atiyah class~\cite{C}. The class is a map $i_*\cO_X\rightarrow i_*\Omega_X[1]$. Let $p_1:X\times X \rightarrow X$ be the projection onto the first factor. Tensoring the map with $p_1^*(E)$ and then pushing forward by $p_{2}$, we get a map $E\rightarrow E\otimes\Omega_X[1]$ which is nothing but the Atiyah class of $E$.

\paragraph\label{Atiyah}
It is easy to see that $\Phi^1_{\varphi}(E)=\mathsf{id}_{\varphi^*E}\otimes\Phi^1_{\varphi}(\cO_X)$ because of the projection formula. The resolution $T_{E,\varphi}$ of $\mu_*E$ has a nice description in terms of the map $\Phi^1_{\varphi}(E)$. The truncation $\tau^{\geq k}T_{E,\varphi}$ of the resolution complex $T_{E,\varphi}$ gives an exact sequence
$$0\rightarrow \mu_* ({N^\vee_{X/S}}^{\otimes k}\otimes E)\rightarrow \varphi^*( {N^\vee_{X/S}}^{\otimes k-1}\otimes E)\rightarrow\cdots\rightarrow \varphi^* ({N^\vee_{X/S}}\otimes E)\rightarrow\varphi^* E\rightarrow\mu_*E \rightarrow 0$$
which defines a map $\Phi^k_{\varphi}(E):\mu_* E \rightarrow\mu_* ({N^\vee_{X/S}}^{\otimes k}\otimes E)[k]$. Due to the construction of the resolution, the map $\Phi^k_{\varphi}(E)$ is equal to the following composite map

$$(\id_{\varphi^* {N^\vee_{X/S}}^{\otimes k-1} }\otimes\Phi^1_{\varphi}(E))\circ(\id_{\varphi^* {N^\vee_{X/S}}^{\otimes k-2} }\otimes\Phi^1_{\varphi}(E))\circ\cdots\circ\Phi^1_{\varphi}(E).$$

Pushing forward the map above onto $S$, we get a map $i_*E\rightarrow i_* ({N^\vee_{X/S}}^{\otimes k}\otimes E)[k]$. Compose it with the natural projection ${N^\vee_{X/S}}^{\otimes k}\rightarrow \wedge^k N^\vee_{X/S}$. We get a map $\Psi^k_{\varphi}(E):i_*E\rightarrow i_*(\wedge^k {N^\vee_{X/S}}\otimes E)[k]$.
\begin{Proposition}\label{prop cal}
In the same setting of Theorem A, consider the isomorphism $I_{\varphi}(E):\mu^*\mu_{*}(E)\cong\mathrm{T}(N^\vee_{X/S}[1])\otimes E$ constructed from a fixed splitting $\varphi$. Due to the adjunction of the functors $\mu^*\dashv\mu_*$, we obtain a map $\Phi_{\varphi}(E): \mu_{*}(E)\rightarrow\mu_*(\mathrm{T}(N^\vee_{X/S}[1])\otimes E)$. Each degree $k$ component of the map $\Phi_{\varphi}(E)$ is the map $\Phi^k_{\varphi}(E)$ defined in Paragraph~\ref{Atiyah}. Similarly, one can construct a map $\Psi_{\varphi}(E):i_*E\rightarrow i_*(\Sym(N^\vee_{X/S}[1])\otimes E)$. Each degree $k$ component of the map $\Psi_{\varphi}(E)$ is $\frac{1}{k!}\Psi^k_{\varphi}(E)$.
\end{Proposition}

\begin{Proof}
The proposition above has been proven~\cite{C} in the case where $i$ is the diagonal embedding, $\varphi$ is $\pi_2$, and $E$ is trivial. The proof in~\cite{C} does not use anything special about the diagonal map. We write the general proof here because the proof will be used throughout this paper.

Consider the unit map $\eta:\mu_*E\rightarrow\mu_*\mu^*\mu_* E$ of the adjunction $\mu^*\dashv\mu_*$. It can be viewed as a map $\eta:\mu_*E=\mu_* E\otimes\cO_{X^{(1)}_{S}}\rightarrow\mu_*\mu^*\mu_* E\cong\mu_*E\otimes\mu_*\cO_X$ where the isomorphism is due to the projection formula. This map is precisely $\id_{\mu_*E}\otimes \varepsilon $, where $\varepsilon:\cO_{X^{(1)}_S}\rightarrow\mu_*\cO_X$ is the natural map of algebras. Due to the adjunction $\mu^*\dashv\mu_*$, we have the equality $\eta\circ\mu_*(I_\varphi)=\Phi_{\varphi}(E)$. The isomorphism $I_{\varphi}$ is defined by identifying the resolution complex $T_{E,\varphi}=(\mathrm{T}^c(\varphi^*N^\vee_{X/S}[1])\otimes\varphi^* E,d)$ with $\mu_*E$ in the derived category of $X$. Under this identification, one can show that the map $\eta\circ\mu_*(I_{\varphi})=\Phi_{\varphi}$
is the map below

$$\mbox{\footnotesize
\xymatrix{
 \ar[r] & \varphi^* ({N^\vee_{X/S}}^{\otimes k}\otimes E)\ar[r]\ar[d] &  \varphi^*( {N^\vee_{X/S}}^{\otimes k-1}\otimes E)\ar[r]\ar[d] & \cdots\ar[r]\ar[d] & \varphi^* ({N^\vee_{X/S}}\otimes E)\ar[r]\ar[d] & \varphi^* E\ar[d]\\
 \ar[r] & \mu_* ({N^\vee_{X/S}}^{\otimes k}\otimes E)\ar[r]^{0} &  \mu_*( {N^\vee_{X/S}}^{\otimes k-1}\otimes E)\ar[r]^{~~~~~~~~0} & \cdots\ar[r]^{0~~~~~~~~} & \mu_* ({N^\vee_{X/S}}\otimes E)\ar[r]^{~~~~~~~~0} & \mu_* E
}
}
$$
which is the natural map between the chain complexes $T_{E,\varphi}$ and the complex $$\mu_*(\mathrm{T}(N^\vee_{X/S}[1])\otimes E)$$
with trivial differential. The map $\Phi_{\varphi}$ above factors through the truncation $\tau^{\geq-k }T_{E,\varphi}$ as follows

$$\mbox{\footnotesize
\xymatrix{
 \ar[r] & \varphi^* ({N^\vee_{X/S}}^{\otimes k}\otimes E)\ar[r]\ar[d] &  \varphi^*( {N^\vee_{X/S}}^{\otimes k-1}\otimes E)\ar[r]\ar[d] & \cdots\ar[r]\ar[d] & \varphi^* ({N^\vee_{X/S}}\otimes E)\ar[r]\ar[d] & \varphi^* E\ar[d]\\
 \ar[r] & \mu_* ({N^\vee_{X/S}}^{\otimes k}\otimes E)\ar[r]\ar[d] &  \varphi^*( {N^\vee_{X/S}}^{\otimes k-1}\otimes E)\ar[r]\ar[d] & \cdots\ar[r]\ar[d] & \varphi^* ({N^\vee_{X/S}}\otimes E)\ar[r]\ar[d] & \varphi^* E\ar[d]\\
 \ar[r] & \mu_* ({N^\vee_{X/S}}^{\otimes k}\otimes E)\ar[r]^{0} &  \mu_*( {N^\vee_{X/S}}^{\otimes k-1}\otimes E)\ar[r]^{~~~~~~~~0} & \cdots\ar[r]^{0~~~~~~~~} & \mu_* ({N^\vee_{X/S}}\otimes E)\ar[r]^{~~~~~~~~0} & \mu_* E
}
}
$$
We look at the degree $k$-th component of the map $\Phi_{\varphi}$, i.e.,
$$\mbox{\footnotesize
\xymatrix{
 \ar[r] & \varphi^* ({N^\vee_{X/S}}^{\otimes k}\otimes E)\ar[r]\ar[d] &  \varphi^*( {N^\vee_{X/S}}^{\otimes k-1}\otimes E)\ar[r]\ar[d] & \cdots\ar[r]\ar[d] & \varphi^* ({N^\vee_{X/S}}\otimes E)\ar[r]\ar[d] & \varphi^* E\ar[d]\\
 \ar[r] & \mu_* ({N^\vee_{X/S}}^{\otimes k}\otimes E)\ar[r]\ar[d] &  \varphi^*( {N^\vee_{X/S}}^{\otimes k-1}\otimes E)\ar[r]\ar[d] & \cdots\ar[r]\ar[d] & \varphi^* ({N^\vee_{X/S}}\otimes E)\ar[r]\ar[d] & \varphi^* E\ar[d]\\
 \ar[r] & \mu_* ({N^\vee_{X/S}}^{\otimes k}\otimes E)\ar[r]^{0} \ar[d] &  \mu_*( {N^\vee_{X/S}}^{\otimes k-1}\otimes E)\ar[r]^{~~~~~~~~0} & \cdots\ar[r]^{0~~~~~~~~} & \mu_* ({N^\vee_{X/S}}\otimes E)\ar[r]^{~~~~~~~~0} & \mu_* E\\
    & \mu_* ({N^\vee_{X/S}}^{\otimes k}\otimes E) .
}
}
$$
The exact sequence on the top is a resolution of $\mu_*E$, so it is identified with $\mu_*E$ in the derived category of $X$. Therefore the composition of the vertical chain maps above is a map $\mu_*E\rightarrow\mu_* ({N^\vee_{X/S}}^{\otimes k}\otimes E)[k]$ which is the degree $k$-th component of $\Phi_\varphi$. One can conclude that the degree $k$-th component of $\Phi_{\varphi}$ is the map $\Phi^k_\varphi$ defined by the truncating the complex $T_{E,\varphi}$ in Paragraph~\ref{Atiyah}. The proof for $\Psi_\varphi$ is similar.
\end{Proof}

\begin{Proposition}\label{prop 2.7}
We pull back the short exact sequence
$$0\rightarrow\mu_*(N^\vee_{X/S}\otimes E)\rightarrow \varphi^*(E)\rightarrow\mu_* E\rightarrow0 $$
by $\mu$.
It defines an exact triangle
$$\mu^*\mu_* (N^\vee_{X/S}\otimes E)\rightarrow E\rightarrow \mu^*\mu_* E\rightarrow \mu^*\mu_*( N^\vee_{X/S}\otimes E)[1].$$
The diagram
$$
\xymatrix{
\mu^*\mu_* E\ar[r] \ar[d]& \mu^*\mu_*( N^\vee_{X/S}\otimes E)[1]\ar[d]\\
E\otimes {N^\vee_{X/S}}^{\otimes k}[k] \ar[r]^{id~~~~~~~~} & E\otimes {N^{\vee}_{X/S}}^{\otimes k-1}[k-1][1]
}
$$
 is commutative, where the vertical maps are the adjunctions to $\Phi_{E}^k$ and $\Phi^{k-1}_{E\otimes N^\vee_{X/S}[1]}$ using the adjunction $\mu^*\dashv\mu_*$ of the functors.

\end{Proposition}

\begin{Proof}
It suffices to prove the commutativity of the diagram
$$
\xymatrix{
\mu_* E\ar[r] \ar[d]^{\Phi_{E}^k}& \mu_*( N^\vee_{X/S}\otimes E)[1]\ar[d]^{\Phi^{k-1}_{E\otimes N^\vee_{X/S}[1]}}\\
\mu_*(E\otimes {N^\vee_{X/S}}^{\otimes k})[k] \ar[r]^{id~~~~~~~~} & \mu_* (E\otimes N^\vee_{X/S}[1]\otimes {N^{\vee}_{X/S}}^{\otimes k-1})[k-1]
}
$$
which follows immediately from the construction of the resolution $T_{E,\varphi}$ of $\mu_* E$. The resolution complex is defined by stringing together a family of short exact sequences.\end{Proof}

The proposition above shows that the study of $\Phi^k_{E,\varphi}$ can be reduced to the study of $\Phi^{0}_{E\otimes {N^\vee_{X/S}}^{\otimes k},\varphi }$. The second one is easier because there is no cohomological shift.

\begin{Proposition}\label{prop 2.8}
There is an isomorphism of
$$
\xymatrix{
E\ar[r]\ar[dd]^{\id} & \mu^*\mu_*(E)\ar[r] \ar[dd]^{I_{\varphi}(E)}_{\cong} & \mu^*\mu_*(E\otimes N^\vee_{X/S}[1]) \ar[dd]^{I_{\varphi}(E\otimes N^\vee_{X/S}[1])}_{\cong} \\
&  & \\
E\ar[r] & \mathrm{T}(N^\vee_{X/S}[1])\otimes E    \ar[r]  & \mathrm{T}(N^\vee_{X/S}[1])\otimes(E\otimes N^\vee_{X/S})[1]
}
$$
exact triangles. The vertical maps are the isomorphisms $I_\varphi$ applied to $E$ and $E\otimes N^\vee_{X/S}[1]$. The exact triangle on the top is defined in Proposition~\ref{prop 2.7}. The kernel of the natural projection map $\mathrm{T}(N^\vee_{X/S}[1])\otimes E\rightarrow \mathrm{T}(N^\vee_{X/S}[1])\otimes(E\otimes N^\vee_{X/S})[1]$ is E. Therefore the bottom line of the diagram above forms a short exact sequence of complexes which can be viewed as an exact triangle in the derived category.
\end{Proposition}
\begin{Proof}
Take the direct sum of the maps $\Phi^k_{E}$ and $\Phi^{k-1}_{E\otimes N^\vee_{X/S}[1]}$ in the proof of Proposition~\ref{prop 2.7} for all $k\geq1$. And notice that the quotient map $\mu^*\mu_* E\rightarrow E$ naturally splits by the map $E\rightarrow\mu^*\mu_* E$ constructed in Proposition~\ref{prop 2.7}.
\end{Proof}

\begin{Proof}[Proof of Theorem B (1).]
For any vector bundle $F$ on $X^{(1)}_S$, tensor it with the short exact sequence
$$0\rightarrow\mu_*N^\vee_{X/S}\rightarrow \cO_{X^{(1)}_S}\rightarrow\mu_*\cO_X\rightarrow0.$$

We get
$$0\rightarrow\mu_*(N^\vee_{X/S}\otimes F|_X)\rightarrow F\rightarrow\mu_*(F|_X)\rightarrow0. $$
Pushing forward by the first order splitting $\varphi$, we get
$$0\rightarrow(N^\vee_{X/S}\otimes F|_X)\rightarrow \varphi_*F\rightarrow(F|_X)\rightarrow0. $$

It is known that the exact sequence splits if and only if $F$ is isomorphic to $\varphi^*(F|_X)$~\cite{CG}.

Choose $F$ and $\varphi$ above to be $\varphi^*_1E$ and $\varphi_2$ respectively. Then the class $\alpha_{\varphi_1,\varphi_2}(E)$ is zero if and only if $\varphi^*_1E$ is isomorphic to $\varphi^*_2E$ because of the reason above.

Denote the two HKR isomorphisms $i^*i_*(E)\cong\Sym(N^\vee_{X/S}[1])\otimes E$ constructed from the two splittings $\varphi_1$ and $\varphi_2$ by $\HKR_{\varphi_1}(E)$ and $\HKR_{\varphi_2}(E)$. The two isomorphism are equal over $X\times X$ is equivalent to $\HKR_{\varphi_1}(E)=\HKR_{\varphi_2}(E)$ for all $E$. They are equal over $X$ is equivalent to $\HKR_{\varphi_1}(\cO_X)=\HKR_{\varphi_2}(\cO_X)$.

We consider the isomorphism over $X\times X$. We prove that the class $\alpha_{\varphi_1,\varphi_2}(E)$ vanishes for all $E$ if $\HKR_{\varphi_1}$ is equal to $\HKR_{\varphi_2}$ over $X\times X$. If $\HKR_{\varphi_1}$ is equal to $\HKR_{\varphi_2}$ over $X\times X$, then the diagram
$$
\xymatrix{
i^*i_*(E)\ar[r]^{\id}\ar[d]^{\HKR_{\varphi_1}} & i^*i_*(E)\ar[d]^{\HKR_{\varphi_2}}\\
\Sym(N^\vee_{X/S}[1])\otimes(E)\ar[r]^{id}\ar[d] & \Sym(N^\vee_{X/S}[1])\otimes(E)\ar[d]\\
N^\vee_{X/S}\otimes E[1]\ar[r]^{\id} & N^{\vee}_{X/S}\otimes E[1]\\
}
$$
is commutative, where the vertical map $\Sym(N^\vee_{X/S}[1])\otimes(E)\rightarrow N^\vee_{X/S}\otimes E[1] $ is the natural projection. Due to Proposition~\ref{prop cal}, we know that the two composite vertical maps $i^*i_*(E)\rightarrow N^\vee_{X/S}\otimes E[1]$ in the diagram above are adjunction to $\Psi^1_{\varphi_1}$ and $\Psi^1_{\varphi_2}$ respectively. The commutativity of the diagram above is equivalent to saying that $\Psi^1_{\varphi_1}$ and $\Psi^1_{\varphi_2}$ are equal, i.e., there is an isomorphism between the two short exact sequence
$$\xymatrix{
0\ar[r] & i_*(N^\vee_{X/S}\otimes E)\ar[r]\ar[d]^{\id} &\nu_* \varphi_1^* E \ar[r]\ar[d] & i_* E\ar[r]\ar[d]^{\id} & 0\\
0\ar[r] & i_*(N^\vee_{X/S}\otimes E)\ar[r] &\nu_* \varphi_2^* E \ar[r]  & i_* E\ar[r] & 0.
}$$
It is enough to consider the exact sequence on $X^{(1)}_S$
$$\xymatrix{
0\ar[r] & \mu_*(N^\vee_{X/S}\otimes E)\ar[r]\ar[d]^{\id} & \varphi_1^* E \ar[r]\ar[d] & \mu_* E\ar[r]\ar[d]^{\id} & 0\\
0\ar[r] & \mu_*(N^\vee_{X/S}\otimes E)\ar[r] & \varphi_2^* E \ar[r]  & \mu_* E\ar[r] & 0
}$$
instead of on $S$. The commutativity of the diagram above implies that $\varphi^*_1E\cong\varphi^*_2E$, equivalently, the class $\alpha_{\varphi_1,\varphi_2}(E)$ is zero for any $E$.

We prove that $\HKR_{\varphi_1}$ is equal to $\HKR_{\varphi_2}$ over $X\times X$ if the class $\alpha_{\varphi_1,\varphi_2}(E)$ vanishes for all $E$. In particular, the class $\alpha_{\varphi_1,\varphi_2}(E\otimes N^\vee_{X/S})$ vanishes, and then there is an isomorphism $\varphi_1^*(N^\vee_{X/S}\otimes E)\cong\varphi^*_2 (N^\vee_{X/S}\otimes E)$ which pull back to identity map of $N^{\vee}_{X/S}\otimes E$ by $\mu$. One sees that the isomorphism induces an isomorphism of complexes between $T_{E,\varphi_1}$ and $T_{E,\varphi_2}$ which pull back to the identity map on $\mathrm{T}(N^\vee_{X/S}[1])\otimes E$ by $\mu$.  Therefore we can conclude that $\HKR_{\varphi_1}$ is equal to $\HKR_{\varphi_2}$ over $X\times X$.

We consider the isomorphism over $X$. From the discussion above, one can conclude that $\HKR_{\varphi_1}$ is equal to $\HKR_{\varphi_2}$ over $X$ if $\alpha_{\varphi_1,\varphi_2}(N^\vee_{X/S})$ vanishes because the two resolutions of $\mu_*\cO_X$ build from the two splittings are isomorphic.
\end{Proof}

\paragraph\label{diagonal} Consider the example when $i$ is the diagonal embedding $X\hookrightarrow X\times X=S$ and $\varphi_i$ are $\pi_i$ for $i=1,2$. Then $\alpha_{\pi_1,\pi_2}(E)$ is the Atiyah class $\mathsf{at}(E)$ of $E$, so we can conclude that $\HKR_{\pi_1}$ is not equal to $\HKR_{\pi_2}$ over $X\times X$ for general $X$.



\begin{Corollary}\label{cor first order over X}
In the same setting of Theorem A, the class $\alpha_{\varphi_1,\varphi_2}(N^\vee_{X/S})$ is zero if and only if the $I_{\varphi_1}(\cO_X)$ is equal to $I_{\varphi_2}(\cO_X)$.
\end{Corollary}
\begin{Proof}
We know that $I_{\varphi_1}(\cO_X)$ is equal to $I_{\varphi_2}(\cO_X)$ if the class $\alpha_{\varphi_1,\varphi_2}(N^\vee_{X/S})$ is zero from the proof of the first part of Theorem B above. Conversely, if $I_{\varphi_1}(\cO_X)$ is equal to $I_{\varphi_2}(\cO_X)$, then $I_{\varphi_1}(N^\vee_{X/S})$ is equal to $I_{\varphi_2}(N^\vee_{X/S})$ due to Propsition~\ref{prop 2.8}. Then $\HKR_{\varphi_1}(N^\vee_{X/S})$ is equal to $\HKR_{\varphi_2}(N^\vee_{X/S})$. From the proof of the first part of Theorem B, we know that the class $\alpha_{\varphi_1,\varphi_2}(N^\vee_{X/S})$ vanishes.
\end{Proof}

The commutative diagram in Proposition~\ref{prop 2.8} is crucial in the proof above. However, we can not get a similar commutative diagram by replacing $\mu$ by $i$ because there is no map $E\rightarrow i^*i_*(E)$ generally.
\paragraph\label{subalgebra}
From the construction of the HKR isomorphism, we know~\cite{AC} that there is a commutative diagram of algebras
$$\xymatrix{
i^*i_*\cO_X\ar[r]\ar[d]^{\HKR_\varphi} & \mu^*\mu_*\cO_X\ar[d]^{I_\varphi}\\
\Sym(N^\vee_{X/S}[1])\ar[r] & \mathrm{T}^c(N^\vee_{X/S}[1]),
}
$$
where $\mathrm{T}^c(N^{\vee}_{X/S}[1])$ is the free tensor coalgebra with the commutative shuffle product. The symmetric algebra is naturally a subalgebra of the tensor coalgebra. Two splittings produce an automorphism $\HKR_{\varphi_1}\circ\HKR_{\varphi_2}^{-1}$ of $\Sym(N^{\vee}_{X/S}[1])$ and an automorphism $I_{\varphi_1}\circ I_{\varphi_2}^{-1}$ of $\mathrm{T}^c(N^{\vee}_{X/S}[1])$ respectively. The HKR isomorphisms $\HKR_{\varphi_1}(\cO_X)$ and $\HKR_{\varphi_2}(\cO_X)$ are equal is equivalent to saying that the automorphism is the identity on the subalgebra. The maps $I_{\varphi_1}(\cO_X)$ and $I_{\varphi_2}(\cO_X)$ are equal is equivalent to saying that the automorphism is the identity on the tensor coalgebra.

\section{Proof of Theorem A and Theorem B (2)}

From two splittings, we obtain two maps $I_{\varphi_1}(E)\circ I^{-1}_{\varphi_2}(E)$ and $\HKR_{\varphi_1}\circ\HKR^{-1}_{\varphi_2}$. The first map defines an automorphism of $\mathrm{T}(N^\vee_{X/S}[1])\otimes E$ and the second defines an automorphism of $\Sym(N^\vee_{X/S}[1])$. We compute the first automorphism and then we use the result to prove Theorem A and the second part of Theorem B.


\begin{Proposition}\label{auto}
The map $I_{\varphi_1}(E)\circ I^{-1}_{\varphi_2}(E)$ defines an automorphism of $\mathrm{T}(N^\vee_{X/S}[1])\otimes E$. Write the automorphism in the form of a matrix below. Then the matrix is unipotent
$$\mbox{\small
\xymatrix{
  & E & E\otimes N^\vee_{X/S}[1] & E\otimes {N^\vee_{X/S}}^{\otimes 2}[2] & \cdots\\
E &  \id &  \alpha_{\varphi_1,\varphi_2}(E)                &                   \cdots        &    \cdots   \\
E\otimes N^\vee_{X/S}[1] & 0  &   \id                  &     \alpha_{\varphi_1,\varphi_2}(E\otimes N^\vee_{X/S})                   &  \cdots     \\
E\otimes {N^\vee_{X/S}}^{\otimes 2}[2] & 0 & 0 & \id &\\
\cdots & & & &,}
}
$$
where the maps on the diagonal are the identity maps. The $(k+1,k+2)$-th entry in this matrix is $\alpha_{\varphi_1,\varphi_2}(E\otimes {N^\vee_{X/S}}^{\otimes k})$. 
\end{Proposition}
\begin{Proof}

 We can apply the isomorphisms $I_{\varphi_1}$ and $I_{\varphi_2}$ to $E\otimes{N^\vee_{X/S}}^{\otimes k} [k]$. We get an automorphism $$I_{\varphi_1}(E\otimes{N^\vee_{X/S}}^{\otimes k} [k])\circ I_{\varphi_2}^{-1}(E\otimes{N^\vee_{X/S}}^{\otimes k} [k])$$
of $\mathrm{T}(N^\vee_{X/S}[1])\otimes (E\otimes{N^\vee_{X/S}}^{\otimes k} [k])$. Proposition~\ref{prop 2.8} shows that the $(p,q)$-th entry in the new matrix of the automorphism above is the $(k+p,k+q)$-th entry of the matrix of the automorphism $I_{\varphi_1}(E)\circ I_{\varphi_2}^{-1}(E)$ by induction on $k$. Therefore, to compute all the entries in the matrix, it suffices to compute the $(1,k)$-th entry of the automorphism matrix.

The map $E \rightarrow E\otimes {N^\vee_{X/S}}^{\otimes k}[k]$ in the $1$-th row and $k+1$-th column of the matrix is defined as follows.
The complexes $T_{E,\varphi_1}$ and $T_{E,\varphi_2}$ are resolutions of $\mu_*E$
$$
\xymatrix{
\cdots \ar[r] & \varphi_1^*(E\otimes N^\vee_{X/S}) \ar[r] & \varphi_1^*(E)\ar[d]\\
 & & \mu_* E\\
\cdots \ar[r] & \varphi_2^*(E\otimes N^\vee_{X/S}) \ar[r] & \varphi_2^*(E).\ar[u]\\
}
$$
Since the resolution $T_{E,\varphi_1}\rightarrow\mu_* E$ is a quasi-isomorphism, it is invertible in the derived category of $X$. Therefore, we get an isomorphism $J_{\varphi_1,\varphi_2}:T_{E,\varphi_2}\rightarrow T_{E,\varphi_1}$ from the complex on the bottom to the complex on the top of the diagram above. There is a natural map from $T_{E,\varphi_1}$ to the truncation $\tau^{\geq -k}T_{E,\varphi_1}$, so we get a map $T_{E,\varphi_2}\rightarrow T_{E,\varphi_1}\rightarrow\tau^{\geq -k}T_{E,\varphi_1} $

$$
\xymatrix{
0\ar[r] & \mu_*(E\otimes {N^\vee_{X/S}}^{\otimes k})\ar[r] & \cdots \ar[r] & \varphi_1^*(E\otimes N^\vee_{X/S}) \ar[r] & \varphi_1^*(E)\ar[d]\\
\ar[r] & \varphi^*_1(E\otimes {N^\vee_{X/S}}^{\otimes k})\ar[r] \ar[u]& \cdots \ar[r]\ar[u] & \varphi_1^*(E\otimes N^\vee_{X/S}) \ar[r]\ar[u] & \varphi_1^*(E)\ar[d]\ar[u]\\
 & & & & \mu_* E\\
\ar[r] & \varphi^*_2(E\otimes {N^\vee_{X/S}}^{\otimes k})\ar[r] &\cdots \ar[r] & \varphi_2^*(E\otimes N^\vee_{X/S}) \ar[r] & \varphi_2^*(E).\ar[u]
}
$$
The truncation naturally maps to $\mu_*(E\otimes {N^\vee_{X/S}}^{\otimes k})[k]$ and $\varphi_2^*(E)$ naturally maps to $T_{E,\varphi_2}$. Therefore we get a map $\varphi_2^*(E)\rightarrow \mu_*(E\otimes {N^\vee_{X/S}}^{\otimes k})[k] $. Pulling back by $\mu$, we get a map $E\rightarrow \mu^*\mu_*(E\otimes {N^\vee_{X/S}}^{\otimes k})[k]$. There is a map $\mu^*\mu_*(E\otimes {N^\vee_{X/S}}^{\otimes k})[k]\rightarrow E\otimes {N^\vee_{X/S}}^{\otimes k}[k]$ due to the adjunction. Compose the two maps together, we get $E\rightarrow E\otimes {N^\vee_{X/S}}^{\otimes k}[k]$. We want to show that the map is equal to the map in the $1$-th row and $k+1$-th column of the matrix. It follows from the following two facts.
\begin{itemize}
\item The map $\mu^*J_{\varphi_1,\varphi_2}:\mu^*\mu_*(E)\rightarrow\mu^*\mu_*(E)$ is exactly the automorphism $I_{\varphi_1}\circ I_{\varphi_2}^{-1}$.
\item For any vector bundle $F$, consider the map $\varphi_2^*F\rightarrow \mu_*F$ by $\mu^*$. Pull the map back to $X$, we get a map $F\rightarrow\mu^*\mu_*F$. Compose it with the natural map $\mu^*\mu_*F\rightarrow F$ defined by adjunction. The composite map is the identity map on $F$. To prove what we need above, we choose $F=(E\otimes {N^\vee_{X/S}}^{\otimes k})[k]$.
\end{itemize}
From the discussion above, one can conclude that the $(1,1)$-th entry in the matrix is the identity map on $E$.

We compute the $(1,2)$-th entry in the matrix.
There are two exact triangles
$$\mu_*(E\otimes N^\vee_{X/S})\rightarrow \varphi_k^*(E)\rightarrow \mu_*(E)\stackrel{\Phi^1_{\varphi_k}}{\rightarrow}\mu_*(E\otimes N^\vee_{X/S})[1] $$
for $k=1,2$.
Consider the composite map $$\beta:\varphi_2^*E\rightarrow \mu_*E\stackrel{\Phi^1_{\varphi_1}}{\rightarrow}\mu_*(E\otimes N^\vee_{X/S})[1].$$
We get a map $E\rightarrow(E\otimes N^\vee_{X/S})[1]$ by adjunction. We know that this map is the $(1,2)$-th entry in the automorphism matrix due to the discussion above. We need to show that it is equal to $\alpha_{\varphi_1,\varphi_2}(E)$. The class $\alpha_{\varphi_1,\varphi_2}$ is defined by pushing forward the map $\Phi_{\varphi_1}^1$ by $\varphi_{2*}$
 $$E=  \varphi_{2*}\mu_* E\stackrel{\varphi_{2*}(\Phi^1_{\varphi_1})}{\rightarrow}  \varphi_{2*}\mu_*(E\otimes N^\vee_{X/S})[1]=(E\otimes N^\vee_{X/S})[1],$$
i.e., we have the equality $\alpha_{\varphi_1,\varphi_2}=\varphi_{2*}(\Phi^1_{\varphi_1})$.
 Let $\eta:E\rightarrow\varphi_{2*}\varphi^{*}_2E$ be the unit of the adjunction $\varphi^{*}_2\dashv\varphi_{2*}$. The composite map
$$\varphi_{2*}(\beta)\circ\eta: E\rightarrow\varphi_{2*}\varphi^{*}_2E\rightarrow \varphi_{2*}\mu_*E=E\stackrel{\varphi_{2*}(\Phi^1_{\varphi_1})}{\rightarrow}  \varphi_{2*}\mu_*(E\otimes N^\vee_{X/S})[1]$$
is the map adjunction to $\beta$ by the property of the unit map of the adjunction $\varphi^{*}_2\dashv\varphi_{2*}$. Since $E\rightarrow\varphi_{2*}\varphi^{*}_2E\rightarrow \varphi_{2*}\mu_*E=E$ is the identity, we get the desired result.

\end{Proof}

\begin{Proof}[Proof of Theorem A and Theorem B (2).] From the proof of Proposition~\ref{auto} one can conclude that the $(p,p+k)$-th entry $$E\otimes {N^\vee_{X/S}}^{\otimes p-1}\rightarrow E\otimes {N^\vee_{X/S}}^{\otimes p+k-1}[k]$$
of the automorphism matrix is $\varphi_{2*}\Phi^k_{\varphi_1}(E\otimes {N^\vee_{X/S}}^{\otimes p-1})$. Because of the projection formula and the definition of $\Phi_{\varphi}^k(E)$, one can show that for any $E$ $$\varphi_{2*}\Phi_{\varphi_1}^k(E)=\varphi_{2*}\Phi_{\varphi_1}^1(E)\otimes\id\circ\cdots\circ\varphi_{2*}\Phi_{\varphi_1}^1(E)\otimes\id\circ\varphi_{2*}\Phi_{\varphi_1}^1(E)$$
which implies that the $(p,p+k)$-th entry in the matrix of the automorphism $I_{\varphi_1}(E)\circ I^{-1}_{\varphi_2}(E)$ is determined by the $(p,p+1)$-th entry.

When $E$ is the structure sheaf $\cO_X$, the class $\alpha_{\varphi_1,\varphi_2}({N^\vee_{X/S}}^{\otimes{p-1}})$ in the $(p,p+1)$-th entry is determined by the class $\alpha_{\varphi_1,\varphi_2}(N^\vee_{X/S})$ in the $(2,3)$-th entry of the matrix due to Proposition~\ref{prop derivation}. Namely, the entries in the matrix of the automorphism $I_{\varphi_1}(\cO_X)\circ I^{-1}_{\varphi_2}(\cO_X)$ is completely determined by the class $\alpha_{\varphi_1,\varphi_2}(N^\vee_{X/S})$. When the class $\alpha_{\varphi_1,\varphi_2}(N^\vee_{X/S})$ vanishes, the automorphism matrix of $I_{\varphi_1}(\cO_X)\circ I^{-1}_{\varphi_2}(\cO_X)$ is the identity matrix. This provides another proof of Corollary~\ref{cor first order over X}.


The diagram in Paragraph~\ref{subalgebra} shows that $\Sym(N^\vee_{X/S}[1])$ is naturally a subalgebra of the tensor coalgebra $\mathrm{T}^c(N^\vee_{X/S}[1])$ with the shuffle product. Therefore we can conclude that the corresponding automorphism matrix of $\Sym(N^\vee_{X/S}[1])$ is also unipotent. The inclusion splits as vector spaces
$$\Sym({N^\vee_{X/S}[1]})\hookrightarrow \mathrm{T}^c(N^\vee_{X/S}[1])\stackrel{\exp}{\rightarrow}\mathrm{T}(N^\vee_{X/S}[1])\rightarrow \Sym({N^\vee_{X/S}[1]}),$$
where the first arrow is the inclusion and the other maps have been explained in Paragraph~\ref{resolution}. Due to the reasons above, one can compute the entries in the automorphism matrix of $\HKR_{\varphi_1}(\cO_X)\circ\HKR_{\varphi_2}^{-1}(\cO_X)$ explicitly. The $(p,p+k)$-th entry
$$  \wedge^{p-1}{N^\vee_{X/S}}\rightarrow  \wedge^{p+k-1}{N^\vee_{X/S}}[k]$$
is the composite map
$$ \wedge^{p-1}{N^\vee_{X/S}}\hookrightarrow  {N^\vee_{X/S}}^{\otimes p-1}$$
$$\xymatrix{ \ar[rrr]^{\varphi_{2*}\Phi^k_{\varphi_1}( {N^\vee_{X/S}}^{\otimes p-1})~~~~~~~~} &  &  &  {N^\vee_{X/S}}^{\otimes k+p-1}[k]
}$$
$$\stackrel{\frac{1}{(k+p-1)!}}{\longrightarrow} {N^\vee_{X/S}}^{\otimes k+p-1}[k]\rightarrow\wedge^{p+k-1}{N^\vee_{X/S}}[k].$$
The $(2,3)$-th entry is the class $\alpha^{\mathsf{antisym}}_{\varphi_1,\varphi_2}(N^\vee_{X/S})$. Similarly, one can show that the automorphism matrix of $\HKR_{\varphi_1}(\cO_X)\circ\HKR^{-1}_{\varphi_2}(\cO_X)$ is the identity matrix if $\alpha^{\mathsf{antisym}}_{\varphi_1,\varphi_2}(N^\vee_{X/S})$ vanishes.
\end{Proof}

When $X$ is of codimension two in $S$, Grivaux~\cite[Theorem 4.17]{Gri'} showed that the matrix of the automorphism $\HKR_{\varphi_1}\circ\HKR^{-1}_{\varphi_2}$ is
$$
\begin{bmatrix}
\id & 0 & 0\\
0 & \id & \theta(\chi) \\
0 & 0 & \id  \\
\end{bmatrix},
$$
where the definition of the class $\theta(\chi)$ can be found in {\em loc. cit.}.  Because of Proposition~\ref{prop universal}, this class $\theta(\chi)$ is exactly $\alpha^{\mathsf{antisym}}_{\pi_1,\pi_2}(N^\vee_{X/S})$ which shows that our computation in the proof of Theorem A agrees with Grivaux's result above.

Grivaux also obtained the following theorem~\cite[Theorem 1.2]{Gri'}. If either the conormal bundle $N^\vee_{X/S}$ carries a global holomorphic connection or the map $\varphi_1-\varphi_2$ is an isomorphism between $\Omega_X$ and $N^\vee_{X/S}$, then $\HKR_{\varphi_1}$ and $\HKR_{\varphi_2}$ are equal. In the first situation, the existence of a holomorphic connection implies that the Atiyah class of $N^\vee_{X/S}$ vanishes, which implies $\alpha^{\mathsf{antisym}}_{\pi_1,\pi_2}(N^\vee_{X/S})$ vanishes due to Proposition~\ref{prop universal}. In the second situation, the conormal bundle $N^{\vee}_{X/S}$ is identified with $\Omega_X$ by the isomorphism $\varphi_1-\varphi_2$. The class $\alpha_{\pi_1,\pi_2}(N^\vee_{X/S})$ is nothing but the Atiyah class of $\Omega_X$ in this case. As a consequence, we know that $\alpha^{\mathsf{antisym}}_{\pi_1,\pi_2}(N^\vee_{X/S})$ vanishes because the Atiyah class of $\Omega_X$ is symmetric. This explains that our Theorem B implies Grivaux's theorem.

\begin{Corollary}
Fix a vector bundle $E$, the isomorphisms $I_{\varphi_1}(E)$ and $I_{\varphi_2}(E)$ are equal if and only if $\alpha_{\varphi_1,\varphi_2}(E)$ and $\alpha_{\varphi_1,\varphi_2}(N^\vee_{X/S})$ vanish.
\end{Corollary}

\begin{Proof}
It is clear from the proof of Theorem A and Proposition~\ref{prop derivation}. The $(p,p+k)$-th entry in the matrix of the automorphism $I_{\varphi_1}(E)\circ I_{\varphi_2}^{-1}(E)$ is determined by the $(p,p+1)$-th entry. The $(p,p+1)$-th entry vanishes for all $p$ if and only if $\alpha_{\varphi_1,\varphi_2}(E)$ and $\alpha_{\varphi_1,\varphi_2}(N^\vee_{X/S})$ vanish.
\end{Proof}

\begin{Corollary}\label{cor over X diagonal}
Let $\Delta:X\hookrightarrow X\times X$ be the diagonal embedding, $\pi_1$ and $\pi_2$ be the two first order splittings defined by the two projections. Then $\HKR_{\pi_1}$ is equal to $\HKR_{\pi_2}$ over $X$.
\end{Corollary}

\begin{Proof}
Recall that $E=\cO_X$ and the conormal bundle $N^\vee_{X/S}$ is the cotangent bundle $\Omega_X$ in this case. The Atiyah class $\mathsf{at}(\Omega_X)=\alpha_{\pi_1,\pi_2}(N^\vee_{X/S})$ is symmetric, i.e., it can be viewed as a map $\Omega_X\rightarrow(\Sym^2\Omega_X)[1]$. In this case the class $\alpha^{\mathsf{antisym}}_{\pi_1,\pi_2}(N^\vee_{X/S})$ in Theorem B always vanishes.
\end{Proof}

\paragraph In the case of diagonal embedding, as mentioned in the introduction, the most widely used HKR isomorphism defined by complete bar resolution is equal to the HKR isomorphism defined by $\pi_1$ and $\pi_2$. Applying $\Hom(-,\cO_X)$ to the isomorphism we get the induced isomorphism of vector spaces
$$\HKR_{\pi_2}:\HT^*(X)=\bigoplus_{p+q=*}H^p(X,\wedge^q T_X)\stackrel{\cong}{\rightarrow}\HH^*(X),$$
where the right hand side is the Hochschild cohomology of $X$. There are natural algebra structures on both sides of the isomorphism above, but the HKR isomorphism is not an isomorphism of algebras. Kontsevich~\cite{K} has modified the HKR isomorphism above to obtain an isomorphism of algebras. He defined an automorphism of $\Td^{-\frac{1}{2}}:\HT^*(X)\rightarrow\HT^*(X)$ given by the contraction with the Todd class of $X$. Then the composite map $\HKR_{\pi_2}\circ\Td^{-\frac{1}{2}}$ is an isomorphism of algebras. We show that this composite map is not equal to the HKR isomorphism defined by any first order splitting.

\begin{Corollary}
In the case of diagonal embedding, the map $\HKR_{\pi_2}\circ\Td^{-\frac{1}{2}}$ defined by Kontsevich is not equal to $\HKR_{\varphi}$ for any first order splitting $\varphi$ in general.
\end{Corollary}

\begin{Proof}
We look at the automorphism $\Td^{-\frac{1}{2}}: \HT^*(X)\rightarrow\HT^*(X)$. For a general $X$, the map
  $$H^p(X,T_X)\hookrightarrow\HT^*(X)=\bigoplus_{p+q=*}H^p(X,\wedge^q T_X)\stackrel{\Td^{-\frac{1}{2}}}{\rightarrow}\HT^*(X)\rightarrow H^{p+1}(X,\cO_X)$$
is nonzero because it is the contraction with the first Chern class of $X$. In particular, for a general $X$, the map is nonzero when $p=0$.

However, for any first order splitting $\varphi$, the $(1,2)$-entry in the automorphism matrix of $\HKR_{\varphi}\circ\HKR^{-1}_{\pi_2}$ is zero because the class $\alpha_{\varphi,\pi_2}(\cO_X)$ vanishes. It implies that $\HKR_{\pi_2}\circ\Td^{-\frac{1}{2}}$ can not be equal to $\HKR_{\varphi}$ for any $\varphi$.
\end{Proof}


\end{document}